\documentclass[12pt]{article}
\usepackage[latin1]{inputenc}
\newtheorem{...}{...}
\begin{document}
\title{Mapping Among Line Elements}
\author{A. C. V. V. de Siqueira}
\date{}
\maketitle
\begin{abstract}
In this paper, we revisit our paper "Matrix Riccati Equations, Kaluza-Klein, Finsler Spaces, and Mapping Among Manifolds' \cite{1}. We will build mapping among generalized quadratic Hamiltonians and we construct Calabi's Riemannian Line Elements for non-quadratic and generalized Hamiltonians. As an application, we use conformally flat forms of two general pseudo-Riemannian line elements embedded in two flat manifolds and obtain an analytical and exact solution of the mapping between these two manifolds as well as an infinite set of exact solutions of the associated matrix Riccati equation.
\end{abstract}
\vspace{3cm}
${}^*$ E-mail: antonio.vsiqueira@ufrpe.br
\newline
${}^*$ E-mail: siqueira.acvv@gmail.com
\newpage
\renewcommand{\theequation}{\thesection.\arabic{equation}}
\section{ Introduction} 
Mapping among manifolds can be useful to solving systems of differential equations.
\newline
This paper is organized as follows, in Sec. $2$ and $3$, we describe some results of mapping presented in a previous paper \cite{1}. In Sec. $4$, we generalize mapping among generalized quadratic Hamiltonians \cite{2}. In Sec. $5$, we show conformally flat forms of a general pseudo-Riemannian line element \cite{1}. In Sec. $6$, we embed an n-dimensional conformally flat form of a pseudo-Riemannian manifold in a flat manifold of dimension $n+2$ \cite{1}. In Sec. $7$, we construct Calabi's Riemannian line elements for non-quadratic and generalized Hamiltonians. In Sec. $8$, we obtain an analytical and exact solution of mapping between two generalized pseudo-Riemannian manifolds. In Sec. $9$, we construct an infinite set of exact solutions of the associated matrix Riccati equation. In Sec. $10$, we summarize the main results of this paper.
\renewcommand{\theequation}{\thesection.\arabic{equation}}
\section{\bf Modified Hamiltonian Formalism}
\setcounter{equation}{0} $         $
In this section we consider some results presented in  \cite{1}.
\newline
Consider a time-dependent Hamiltonian
$H({\tau})$ in which ${\tau}$ is an affine parameter, proper-time, for example. 
\newline
Let us define 2n variables, which will be called ${\xi}^j$, with index j running from 1 to 2n so that
we have
${\xi}^j$ $\in$
$({\xi}^1,\ldots,{\xi}^n,{\xi}^{n+1},\ldots,{\xi}^{2n})$ =$(
{q}^1,\ldots,{q}^n,{p}^1,\ldots,{p}^n)$ in which ${q}^j$ and ${p}^j$
are coordinates and momenta, respectively. 
\newline
We now define the Hamiltonian by
\begin{equation}
 H({\tau})=\frac{1}{2}H_{ij}{\xi}^i{\xi}^j,
\end{equation}
in which $H_{ij}$ is a  symmetric matrix. 
\newline
We impose that the Hamiltonian obeys the Hamilton equation
\begin{equation}
\frac{d{\xi}^i}{d\tau }={J}^{ik}\frac{\partial{H}}{\partial{\xi}^k}.
\end{equation}
Equation (2.2) introduces the symplectic J given by
\begin{equation}
\left(%
\begin{array}{cc}
  O & I \\
  -I & O \\
\end{array}%
\right)
\end{equation}
in which O and I are the $n \textbf{x}n$ zero and identity matrices, respectively. 
\newline
We now make a linear transformation from ${\xi}^j$ to ${\eta}^j$ given by
\begin{equation}
{\eta}^j={{T}^j}_k{\xi}^k,
\end{equation}
in which ${{T}^j}_k$ could be a non-symplectic matrix and the new Hamiltonian is given by
\begin{equation}
 Q=\frac{1}{2}Q_{ij}{\eta}^i{\eta}^j,
\end{equation}
in which $Q_{ij}$ is a symmetric matrix. 
\newline
The matrices $H_{ij}$, $Q_{ij}$, and ${{T}^j}_k$ obey the following system
\begin{equation}
\frac{d{{T}^i}_j}{d\tau}+{{T}^i}_k{J}^{kl}X_{lj}=J^{im}\frac{d{t}}{d\tau}Y_{ml}{{T}^j}_k,
\end{equation}
in which $2X_{lj}=\frac{\partial{H_{ij}}}{\partial{\xi}^l
}\xi^{i}+2H_{lj}$, $2Y_{ml}=\frac{\partial{Q_{il}}}{\partial{\eta}^m
}\eta^{i}+2Q_{ml},$ t and $\tau$ can be the proper-times of the particle in two different manifolds. 
\newline
Consider $X_{lj}=Z_{lj}$ and $\frac{d{t}}{d\tau}Y_{ml}=\overline{Y}_{ml}$. Then, (2.6) can be rewritten in the following matrix form
\begin{equation}
 \frac{d{T}}{d\tau}+TJZ=J\overline{Y}T,
\end{equation}
in which $T$, $Z$ and $\overline{Y}$ are  $2n \textbf{x}2n$ matrices as
\begin{equation}
\left(%
\begin{array}{cc}
  T_{1} & T_{2} \\
  T_{3} & T_{4} \\
\end{array}%
\right)
\end{equation}
with similar expressions for $Z$ and $\overline{Y}$. 
\newline
Let us  write (2.7) as follows
\begin{equation}
 \dot{T_1}=\overline{Y}_{3}T_{1}+\overline{Y}_{4}T_{3}+T_{2}Z_{1}-T_{1}Z_{3},
\end{equation}
\begin{equation}
 \dot{T_2}=\overline{Y}_{3}T_{2}+\overline{Y}_{4}T_{4}+T_{2}Z_{2}-T_{1}Z_{4},
\end{equation}
\begin{equation}
 \dot{T_3}=-\overline{Y}_{1}T_{1}-\overline{Y}_{2}T_{3}+T_{4}Z_{1}-T_{3}Z_{3},
\end{equation}
\begin{equation}
 \dot{T_4}=-\overline{Y}_{1}T_{2}-\overline{Y}_{2}T_{4}+T_{4}Z_{2}-T_{3}Z_{4}.
\end{equation}
Now consider
\begin{equation}
 \dot{S_1}=\overline{Y}_{3}S_{1}+\overline{Y}_{4}S_{3},
\end{equation}
\begin{equation}
 \dot{S_2}=\overline{Y}_{3}S_{2}+\overline{Y}_{4}S_{4},
\end{equation}
\begin{equation}
 \dot{S_3}=-\overline{Y}_{1}S_{1}-\overline{Y}_{2}S_{3},
\end{equation}
\begin{equation}
 \dot{S_4}=-\overline{Y}_{1}S_{2}-\overline{Y}_{2}S_{4},
\end{equation}
and
\begin{equation}
\dot{R_1}=R_{2}Z_{1}-R_{1}Z_{3},
\end{equation}
\begin{equation}
\dot{R_2}=R_{2}Z_{2}-R_{1}Z_{4},
\end{equation}
\begin{equation}
\dot{R_3}=R_{4}Z_{1}-R_{3}Z_{3},
\end{equation}
\begin{equation}
\dot{R_4}=R_{4}Z_{2}-R_{3}Z_{4}.
\end{equation}
The systems (2.9)-(2.12),(2.13)-(2.16), and(2.17)-(2.20) can be placed on a compact form as follows
\begin{equation}
\dot{S}=J{\bar{Y}}S
\end{equation}
\begin{equation}
\dot{R}=-RJZ,
\end{equation}
\begin{equation}
{T}={SAR},
\end{equation}
In which matrix ${A}$ is constant and $2nX2n$. 
\newline
A more explicit form for (2.23) is given by
\begin{equation}
 {T_1}=(S_{1}a+S_{2}b)R_{1}+(S_{1}d+S_{2}c)R_{3},
\end{equation}
\begin{equation}
{T_2}=(S_{1}a+S_{2}b)R_{2}+(S_{1}d+S_{2}c)R_{4},
\end{equation}
\begin{equation}
{T_3}=(S_{3}a+S_{4}b)R_{1}+(S_{3}d+S_{4}c)R_{3},
\end{equation}
\begin{equation}
{T_4}=(S_{3}a+S_{4}b)R_{2}+(S_{3}d+S_{4}c)R_{4},
\end{equation}
in which a, b, c, and d are constant $n \textbf{x}n$ matrices.
\newline
From the theory of first-order differential equation systems \cite{3}, it is well known that each systems (2.15)-(2.20) have  solutions in the region, in which $Z_{lj}$ and ${\bar{Y}_{ml}}$ are continuous functions.
\newpage
\renewcommand{\theequation}{\thesection.\arabic{equation}}
\section{\bf Generalized Mapping Among Manifolds  }
\setcounter{equation}{0} $         $
In this section, we present more results from \cite{1}, in which the  matrix Riccati equation was introduced.
\newline
We consider a time-dependent function ${H}$
in which ${\tau}$ is an affine parameter. 
\newline
Let us define 2n variables, which will be called ${\xi}^j$, with index j running from 1 to 2n so that we have ${\xi}^j$ $\in$
$({\xi}^1,\ldots,{\xi}^n,{\xi}^{n+1},\ldots,{\xi}^{2n})$ =$(
{q}^1,\ldots,{q}^n,{p}^1,\ldots,{p}^n)$, in which ${q}^j$ and ${p}^j$
may or may not be the usual coordinates and momenta, respectively. 
\newline
We now define the function by
\begin{equation}
 H({\tau})=\frac{1}{2}H_{ij}{\xi}^i{\xi}^j,
\end{equation}
in which $H_{ij}$ is a  symmetric matrix. 
\newline
Consider the following system
\begin{equation}
\frac{d{\xi}^i}{d\tau
}=C^{ik}\frac{\partial{H}}{\partial{\xi}^k } .
\end{equation}
The equation (3.2) introduces the  $2n\textbf{x}2n$ matrix $C$ given by
\begin{equation}
\left(%
\begin{array}{cc}
  C_1 & C_2 \\
  C_3 & C_4 \\
\end{array}%
\right)
\end{equation}
in which the $C_i$ are $n \textbf{x}n$ matrices, which can be functions of usual coordinates and momenta.
\newline
We now make a linear transformation from ${\xi}^j$ to ${\eta}^j$ given by
\begin{equation}
{\eta}^j={{T}^j}_k{\xi}^k,
\end{equation}
in which ${{T}^j}_k$ could be a non-sympletic matrix and the new function
is given by
\begin{equation}
Q=\frac{1}{2}Q_{ij}{\eta}^i{\eta}^j,
\end{equation}
in which $Q_{ij}$ is a  symmetric matrix. 
\newline
Let us consider that (3.5) obeys the following equation
\begin{equation}
\frac{d{\eta}^i}{dt}=B^{ik}\frac{\partial{Q}}{\partial{\eta}^k },
\end{equation}
in which $B$ is given by
\begin{equation}
\left(%
\begin{array}{cc}
  B_1 & B_2 \\
  B_3 & B_4 \\
\end{array}%
\right)
\end{equation}
in which the $B_i$ are $n \textbf{x}n$ matrices, which can be functions of usual coordinates and momenta.
\newline
The matrices  $H_{ij}$,  $Q_{ij}$, and ${{T}^j}_k$ obey the following system
\begin{equation}
 \frac{d{{T}^i}_j}{d\tau}+{{T}^i}_kC^{kl}X_{lj}=B^{im}(\frac{d{t}}{d\tau}Y_{ml}){{T}^j}_k,
\end{equation}
in which  $2X_{lj}=\frac{\partial{H_{ij}}}{\partial{\xi}^l
}\xi^{i}+2H_{lj}$ and
$2Y_{ml}=\frac{\partial{Q_{il}}}{\partial{\eta}^m
}\eta^{i}+2Q_{ml},$ and $t$ and $\tau$ can be the proper-times of the
particle in two different manifolds. 
\newline
Consider $X_{lj}=Z_{lj}$ and $\frac{d{t}}{d\tau}Y_{lm}=\overline{Y}_{lm}$. Then (3.8) can be rewritten in the following matrix form
\begin{equation}
 \frac{d{T}}{d\tau}+TCZ=B\overline{Y}T,
\end{equation}
in which $T$ are  $2n \textbf{x}2n$ matrices as
\begin{equation}
\left(%
\begin{array}{cc}
  T_{1} & T_{2} \\
  T_{3} & T_{4} \\
\end{array}%
\right)
\end{equation}
with similar expressions for $Z$ and $\overline{Y}$. 
\newline
Let us consider the matrices, T, S, A, R, in which each matrix can be functions of coordinates and momenta,
\begin{equation}
\left(%
\begin{array}{cc}
  T_{1} & T_{2} \\
  T_{3} & T_{4} \\
\end{array}%
\right)
\end{equation}
\begin{equation}
\left(%
\begin{array}{cc}
  S_{1} & S_{2} \\
  S_{3} & S_{4} \\
\end{array}%
\right)
\end{equation}
\begin{equation}
\left(%
\begin{array}{cc}
  A_{1} & A_{2} \\
  A_{3} & A_{4} \\
\end{array}%
\right)
\end{equation}
\begin{equation}
\left(%
\begin{array}{cc}
  R_{1} & R_{2} \\
  R_{3} & R_{4} \\
\end{array}%
\right)
\end{equation}
Now consider the following matricial equation,
\begin{equation}
T=SAR.
\end{equation}
The derivative of (3.15) is given by
\begin{equation}
 \dot{T}= \dot{S}AR+T\dot{A}R+TS\dot{R}
\end{equation}
Let us define two systems of matrix Riccati Equations \cite{2},
\begin{equation}
 \dot{S}=B\overline{Y}S+D+SAF,
\end{equation}
\begin{equation}
 \dot{R}=-RCZ+E+GAR.
\end{equation}
For the particular case in which
\begin{equation}
D+SAF=0,
\end{equation}
\begin{equation}
E+GAR=0,
\end{equation} 
\begin{equation}
B=C=J,
\end{equation}
and A is a constant matrix, the systems (3.17)-(3.18) are reduced to the systems (2.21)-(2.22).
\newline
Replacing (3.17) and (3.18) in (3.16) and assuming that
\begin{equation}
 S\dot{A}R+DAR+SAE+S[A(F+G)A]R=0,
\end{equation}
we have the following simplification of (3.16), given by
\begin{equation}
\frac{d{T}}{d\tau}+TCZ-B\overline{Y}T=0.
\end{equation} 
Equation $(3.9)$ is the motion equation, which was obtained from the equations of motion and is identical to $(3.23)$.
\newline
From the theory of first-order differential equation systems \cite{3}, it is well known that system (3.23) has a solution in the region in which $Z_{lj}$ and $\overline{Y}_{ml}$ are continuous functions. 
If S and R are non-singular matrices, we can multiply (3.22) by $S^{-1}$ on the left side and by $R^{-1}$ on the right side, obtaining the matrix Riccati 
\begin{equation}
\dot{A}+S^{-1}DA+AER^{-1}+A(F+G)A=0.
\end{equation}
Notice that the transformation (3.4) is directly associated with matrix Riccati differential equations (3.23), (3.17), (3.18), and (3.24).
\renewcommand{\theequation}{\thesection.\arabic{equation}}
\section{\bf Generalized Quadratic Hamiltonians  }
\setcounter{equation}{0} $         $
In this section, we generalize important results obtained by Leach in \cite{2}, in which mapping is constructed among generalized quadratic Hamiltonians. 
\newline
As in Section $3$, we consider a time-dependent function $H({\tau})$
in which ${\tau}$ is an affine parameter. 
\newline
Let us define 2n variables, which will be called ${\xi}^j$, with index j running from 1 to 2n so that we have ${\xi}^j$ $\in$
$({\xi}^1,\ldots,{\xi}^n,{\xi}^{n+1},\ldots,{\xi}^{2n})$ =$(
{q}^1,\ldots,{q}^n,{p}^1,\ldots,{p}^n)$, in which ${q}^j$ and ${p}^j$
may or may not be the usual coordinates and momenta, respectively. 
\newline
We now define the function by
\begin{equation}
 H({\tau})=\frac{1}{2}H_{ij}({\tau}){\xi}^i{\xi}^j+G_{i}({\tau}){\xi}^i+D({\tau}),
\end{equation}
in which $H_{ij}({\tau})$ is a  symmetric matrix, $G_{i}({\tau})$ are components of a vector, and $D({\tau})$ is a scalar function.
\newline
Consider the following system
\begin{equation}
\frac{d{\xi}^i}{d\tau
}=C^{ik}\frac{\partial{H}}{\partial{\xi}^k } .
\end{equation}
Equation (4.2) introduces the  $2n\textbf{x}2n$ matrix $C$ given by
\begin{equation}
\left(%
\begin{array}{cc}
  C_1 & C_2 \\
  C_3 & C_4 \\
\end{array}%
\right)
\end{equation}
in which the ${C_i}$ are $ n \textbf{x}n $ matrices, which can be functions of usual coordinates and momenta. 
\newline
We now make a linear transformation from ${\xi}^j$ to
${\eta}^j$ given by
\begin{equation}
 {\eta}^j={{T}^j}_k{\xi}^k+{r}^i,
\end{equation}
in which ${r}^i$ is a vector, ${{T}^j}_k$ can be a non-sympletic matrix and the new function
is given by
\begin{equation}
 Q({t})=\frac{1}{2}Q_{ij}({t}){\eta}^i{\eta}^j+E_{i}({t}){\eta}^i+F({t}),
\end{equation}
in which $Q_{ij}(t)$ is a symmetric matrix, $E_{i}(t)$ are components of a vector and $F(t)$ is a scalar function.
\newline
Let us consider that (4.5) obeys the following equation
\begin{equation}
\frac{d{\eta}^i}{dt}=K^{ik}\frac{\partial{Q}}{\partial{\eta}^k },
\end{equation}
in which $K$ is given by
\begin{equation}
\left(%
\begin{array}{cc}
  K_1 & K_2 \\
  K_3 & K_4 \\
\end{array}%
\right)
\end{equation}
and $K_i$ are $n \textbf{x}n$ matrices, which can be functions of usual coordinates and momenta.
\newline
The matrices  $H_{ij}$,  $Q_{ij}$, and ${{T}^j}_k$ obey the following systems
\begin{equation}
 \frac{d{{T}^i}_j}{d\tau}+{{T}^i}_k{C}^{kl}{\bar{Z}}_{lj}={K}^{im}(\frac{d{t}}{d\tau}){\bar{Y}}_{ml}{{T}^j}_k,
\end{equation}
\begin{equation}
 \frac{d{{r}^i}}{d\tau}={K}^{im}\frac{d{t}}{d\tau}({\bar{Y}}_{ml}+{\bar{E}}_{m})-{{T}^i}_k{C}^{kl}{\bar{G}}_{l}.
\end{equation}
The matrices ${\bar{Z}}_{lj}$, ${\bar{Y}}_{ml}$, and the vectors ${\bar{E}}_{m}$, ${\bar{G}}_{l}$  obey the following equations
\begin{equation}
 2{\bar{Z}}_{lj}=\frac{\partial{H_{ij}}}{\partial{\xi}^l
}\xi^{i}+2H_{lj}+2\frac{\partial{G_{j}}}{\partial{\xi}^l
},
\end{equation}
\begin{equation}
2{\bar{Y}}_{ml}=\frac{\partial{Q_{il}}}{\partial{\eta}^m
}\eta^{i}+2Q_{ml}+2\frac{\partial{E_{j}}}{\partial{\eta}^l}, 
\end{equation}
\begin{equation}
{\bar{E}}_{m}={E}_{m}+\frac{\partial{F}}{\partial{\eta}^l},
\end{equation}
\begin{equation}
{\bar{G}}_{j}={G}_{j}+\frac{\partial{D}}{\partial{\xi}^j}.
\end{equation}
\newline
Let us consider the matrices, T, S, A, R, in which each matrix can be functions of coordinates and momenta,
\begin{equation}
\left(%
\begin{array}{cc}
  T_{1} & T_{2} \\
  T_{3} & T_{4} \\
\end{array}%
\right)
\end{equation}
\begin{equation}
\left(%
\begin{array}{cc}
  S_{1} & S_{2} \\
  S_{3} & S_{4} \\
\end{array}%
\right)
\end{equation}
\begin{equation}
\left(%
\begin{array}{cc}
  A_{1} & A_{2} \\
  A_{3} & A_{4} \\
\end{array}%
\right)
\end{equation}
\begin{equation}
\left(%
\begin{array}{cc}
  R_{1} & R_{2} \\
  R_{3} & R_{4} \\
\end{array}%
\right)
\end{equation}
Now consider the following matrix equation,
\begin{equation}
T=SAR.
\end{equation}
The derivative of (4.18) is given by
\begin{equation}
 \dot{T}= \dot{S}AR+T\dot{A}R+TS\dot{R}.
\end{equation}
As in Section $3$, let us define two systems of matrix Riccati Equations \cite{3},
\begin{equation}
 \dot{S}=B\overline{Y}S+D+SAF,
\end{equation}
\begin{equation}
 \dot{R}=-RCZ+E+GAR.
\end{equation}
Replacing (4.20) and (4.21) in (4.19), and assuming that
\begin{equation}
 S\dot{A}R+DAR+SAE+S[A(F+G)A]R=0,
\end{equation}
we have the following simplification of (4.19) given by
\begin{equation}
\frac{d{T}}{d\tau}+TC{\bar{Z}}-K{\bar{Y}}T=0.
\end{equation} 
Equation (3.9) is a matrix Riccati equation obtained from the equations of motion and is identical to (4.23).
\newline
The main difference between the treatment used in Section 3 and this one is the transformation law (3.4) to (4.4), and the presence of equation (4.9). 
\newline
With definitions (4.10) and (4.11) we can use the same arguments as in section (3) for the case of generalized quadratic Hamiltonians.
\newpage
\renewcommand{\theequation}{\thesection.\arabic{equation}}
\section{\bf Calabi's Line Elements }
\setcounter{equation}{0} $         $
In this section we consider an important result obtained by Calabi,    \cite{4}. We use a development  
given in \cite{5}.
\newline
For one convex surface $u(x^i)$ Calabi defined a metric by
\begin{equation}
G_{ij}=\frac{{\partial}^2{u}}{\partial{x}^k\partial{x}^j, }
\end{equation}
with line element given by
\begin{equation}
ds^{2}=G_{ij}dx^{i}dx^{j}.
\end{equation}
We call (5.2) Calabi's Line Element.
\newline
There is a class of Lagrangians and Hamiltonians that are conservative and behave like surfaces. However, we can generally associate the positive power of a Hamiltonian H with the Calabi's line element (5.2) as follows
\begin{equation}
G_{ij}=\frac{{\partial}^2{H^n}}{\partial{x}^k\partial{x}^j },
\end{equation}
with $n=1,2,3....$.
\newline
When H is on the Hamilton-Jacobi form, the line element (5.2) could be flat.
\newline
From (5.2), we have Riemann's and Ricci's tensors,
\begin{equation}
R^{\alpha}{}_{\mu \sigma \nu }=\partial_{\nu} \Gamma_{\mu \sigma
}^{\alpha}-\partial_{\sigma}\Gamma_{\mu \nu}^{\alpha}+\Gamma_{\mu
\sigma}^{\eta}\Gamma_{n \nu}^{\alpha}-\Gamma_{\mu \nu}^{\eta}\Gamma
_{\sigma \eta}^{\alpha},
\end{equation}
\begin{equation}
R_{\mu \nu}=R^{\alpha}{_{\mu \alpha \nu }},
\end{equation}
\newline
It is ease to show that
\begin{equation}
R_{hijk}=G^{lm}({\Gamma}_{ijm}{\Gamma}_{hkl}-{\Gamma}_{ikm}{\Gamma}_{hjl}),
\end{equation}
in which ${\Gamma}_{ijm}$ is the Christoffel symbol of the first kind.
\newline
Calabi developed a Riemannian Geometry.
\newpage
We are only interested in (5.2) and in (5.3) because, from (5.2), we have the following Lagrangian and Hamiltonian, 
\begin{equation}
L=G_{ij}\frac{d{{x}^i}}{ds}\frac{d{{x}^j}}{ds},
\end{equation}
and
\begin{equation}
H=G^{ij}\frac{d{{p}_i}}{ds}\frac{d{{p}_j}}{ds},
\end{equation}
in which $L=H$.
\newline
We can use Calabi's line elements to transform non-quadratic Hamiltonians on quadratic Hamiltonians as (5.8). Information based on second derivatives can be important for many non-quadratic Hamiltonians, in which the use of the Calabi's line elements may be justified.
\renewcommand{\theequation}{\thesection.\arabic{equation}}
\section{\bf Conformally Flat Forms of a Line Elements }
\setcounter{equation}{0} $         $
In this section, we decribe the conformally flat form of line elements in local coordinates presented in \cite{1}.
\newline
Let us write a metric and its associated line element as follows
\begin{equation}
G_{\Lambda\Pi}=E_{\Lambda}^{(\mathbf{A})}E_{\Pi}^{(\mathbf{B})}\eta_{(\mathbf{A})(\mathbf{B})},
\end{equation}
and
\begin{equation}
 ds^2= G_{\Lambda\Pi}du^{\Lambda}du^{\Pi},
\end{equation}
in which $ \eta_{(\mathbf{A})(\mathbf{B})}$ and $
E_{\Lambda}^{(\mathbf{A})}$ are respectively flat metric and vielbein
components.
\newline
We choose each
 $ \eta_{(\mathbf{A})(\mathbf{B})}$ as plus or minus Kronecker's delta function.
\newline
Let us define
\begin{equation}
\overline{Z}^{(\mathbf{A})}= E_{\Lambda}^{(\mathbf{A})}u^{\Lambda}.
\end{equation}
From (6.3), we have
\begin{equation}
u^{\Lambda}=E^{\Lambda}_{(\mathbf{A})}\overline{Z}^{(\mathbf{A})},
\end{equation}
with
\begin{equation}
du^{\Lambda}= E^{\Lambda}_{(\mathbf{A})}d\overline{Z}^{(\mathbf{A})}+dE^{\Lambda}_{(\mathbf{A})}\overline{Z}^{(\mathbf{A})}.
\end{equation}
We can write (6.5) in a compact form
\begin{equation}
du^{\Lambda}= d(E^{\Lambda}_{(\mathbf{A})}\overline{Z}^{(\mathbf{A})}).
\end{equation}
Substituting (6.6) in (6.2), we have 
\begin{eqnarray}
ds^2=G_{\Lambda\Pi}du^{\Lambda}du^{\Pi}=\eta_{(\mathbf{A})(\mathbf{B})}d\overline{Z}^{(\mathbf{A})}d^{(\mathbf{B})}\\
\nonumber+G_{\Lambda\Pi}[\overline{Z}^{(\mathbf{A})}\overline{Z}^{(\mathbf{B})}d(E^{\Lambda}_{(\mathbf{A})})d(E^{\Pi}_{(\mathbf{B})})\\\nonumber+2\overline{Z}^{(\mathbf{A})}E^{\Pi}_{(\mathbf{B})}d\overline{Z}^{(\mathbf{B})}d(E^{\Lambda}_{(\mathbf{A})})].
\end{eqnarray}
Let us put (6.7) in the following form
\begin{eqnarray}
 ds^2=\eta_{(\mathbf{A})(\mathbf{B})}d\overline{Z}^{(\mathbf{A})}d\overline{Z}^{(\mathbf{B})}+G_{\Lambda\Pi}[\overline{Z}^{(\mathbf{A})}\overline{Z}^{(\mathbf{B})}\frac{d(E^{\Lambda}_{(\mathbf{A})})}{ds}\frac{d(E^{\Pi}_{(\mathbf{B})})}{ds}\\
 \nonumber+2\overline{Z}^{(\mathbf{A})}E^{\Pi}_{(\mathbf{B})}\frac{dz^{(\mathbf{B})}}{ds}\frac{d(tE^{\Lambda}_{(\mathbf{A})})}{ds}]ds^2.
\end{eqnarray}
From a simple calculation, we have
\begin{eqnarray}
\eta_{(\mathbf{A})(\mathbf{B})}d\overline{Z}^{(\mathbf{A})}d\overline{Z}^{(\mathbf{B})}=[1-G_{\Lambda\Pi}[\overline{Z}^{(\mathbf{A})}\overline{Z}^{(\mathbf{B})}\frac{d(E^{\Lambda}_{(\mathbf{A})})}{ds}\frac{d(E^{\Pi}_{(\mathbf{B})})}{ds}\\
 \nonumber+2\overline{Z}^{(\mathbf{A})}E^{\Pi}_{(\mathbf{B})}\frac{d\overline{Z}^{(\mathbf{B})}}{ds}\frac{d(E^{\Lambda}_{(\mathbf{A})})}{ds}]]ds^2.
\end{eqnarray}
We now define the function
\begin{eqnarray}
\nonumber \exp(-2\sigma)=[1-G_{\Lambda\Pi}[\overline{Z}^{(\mathbf{A})}\overline{Z}^{(\mathbf{B})}\frac{d(E^{\Lambda}_{(\mathbf{A})})}{ds}\frac{d(E^{\Pi}_{(\mathbf{B})})}{ds}\\
 \nonumber+2\overline{Z}^{(\mathbf{A})}E^{\Pi}_{(\mathbf{B})}\frac{d\overline{Z}^{(\mathbf{B})}}{ds}\frac{d(E^{\Lambda}_{(\mathbf{A})})}{ds}]].
\end{eqnarray}
Multiplying (6.11) by $\exp(2\sigma)$, we obtain
\begin{eqnarray}
ds^2= G_{\Lambda\Pi}du^{\Lambda}du^{\Pi}=\exp(2\sigma) \eta_{(\mathbf{A})(\mathbf{B})}d\overline{Z}^{(\mathbf{A})}d\overline{Z}^{(\mathbf{B})}.
\end{eqnarray}
Let us define
\begin{equation}
\exp(2\Phi)=U^{-2},
\end{equation} 
in which 
\begin{equation}
U=[1+\frac{1}{4}K\eta_{(\mathbf{A})(\mathbf{B})}\bar{Z}^{(\mathbf{A})}\bar{Z}^{(\mathbf{B})}],
\end{equation}
Multiplying (6.11) by (6.14) and $\exp(-2\sigma)$, we have
\begin{eqnarray}
\exp(2\Phi)\eta_{(\mathbf{A})(\mathbf{B})}d\overline{Z}^{(\mathbf{A})}d\overline{Z}^{(\mathbf{B})}=\exp(-2\sigma)G_{\Lambda\Pi}du^{\Lambda}du^{\Pi}
\end{eqnarray}
We conclude that (6.12) is the line element of a pseudo-Riemannian metric in a conformally flat form and (6.15) is a line element of a pseudo-Riemannian metric of constant curvature as a function of the metric $G_{\Lambda\Pi}$.
\section{ Embedding a Conformally Flat Manifold in Flat Manifolds } 
\setcounter{equation}{0} $         $
In this section, we consider the embedding of (6.10) using a procedure also presented in \cite{1}.
\newline
Let us rewrite (6.10)as follows
\begin{eqnarray}
ds^{2}= G_{\Lambda\Pi}du^{\Lambda}du^{\Pi}=\exp(2\sigma) \eta_{(\mathbf{A})(\mathbf{B})}d\overline{Z}^{(\mathbf{A})}d\overline{Z}^{(\mathbf{B})}.
\end{eqnarray}
Defining the transformation of coordinates by,
\begin{equation}
y^{(\mathbf{A})}=\exp(\sigma)\overline{Z}^{(\mathbf{A})},
\end{equation}
with $ (A)=(1,2,3,....,n),$
\begin{equation}
y^{n+1}=\exp(\sigma)(\eta_{(\mathbf{A})(\mathbf{B})}\overline{Z}^{(\mathbf{A})}\overline{Z}^{(\mathbf{B})}-\frac{1}{4}),
\end{equation}
and
\begin{equation}
y^{n+2}=\exp(\sigma)(\eta_{(\mathbf{A})(\mathbf{B})}\overline{Z}^{(\mathbf{A})}\overline{Z}^{(\mathbf{B})}+\frac{1}{4}).
\end{equation}
But 
\begin{equation}
\eta_{(\mathbf{A})(\mathbf{B})}\overline{Z}^{(\mathbf{A})}\overline{Z}^{(\mathbf{B})}=G_{\Lambda\Pi}u^{\Lambda}u^{\Pi}.
\end{equation}
Using (6.3) and (7.5) in (7.2), (7.3) and (7.4),
\begin{equation}
y^{(\mathbf{A})}=\exp(\sigma)E_{\Lambda}^{(\mathbf{A})}u^{\Lambda},
\end{equation}
with $ (A)=(1,2,3,....,n),$
\begin{equation}
y^{n+1}=\exp(\sigma)(G_{\Lambda\Pi}u^{\Lambda}u^{\Pi}-\frac{1}{4}),
\end{equation}
and
\begin{equation}
y^{n+2}=\exp(\sigma)(G_{\Lambda\Pi}u^{\Lambda}u^{\Pi}+\frac{1}{4}).
\end{equation}
It is easy to see that
\begin{equation}
\eta_{\mathbf{A}\mathbf{B}}y^{\mathbf{A}}y^{\mathbf{B}}=0,
\end{equation}
in which
\begin{equation}
\eta_{\mathbf{A}\mathbf{B}}=(\eta_{(\mathbf{A})(\mathbf{B})},\eta_{\mathbf{(n+1),}\mathbf{(n+1)}},\eta_{\mathbf{(n+2),}\mathbf{(n+2)}}),
\end{equation}
with
\begin{equation}
\eta_{\mathbf{(n+1),}\mathbf{(n+1)}}=1,
\end{equation}
and
\begin{equation}
\eta_{\mathbf{(n+2),}\mathbf{(n+2)}=-1}.
\end{equation}
From a simple calculation, we can verify that the line elements are given by
\begin{equation}
ds^2=G_{\Lambda\Pi}du^{\Lambda}du^{\Pi}=\exp(2\sigma)\eta_{(\mathbf{A})(\mathbf{B})}dz^{(\mathbf{A})}dz^{(\mathbf{B})}=
\eta_{\mathbf{A}\mathbf{B}}dy^{\mathbf{A}}dy^{\mathbf{B}}.
\end{equation}
From (7.13), we see that an n-dimensional manifold in local coordinates can be put in a conformally flat form and embedded in an (n + 2)-dimensional flat manifold.
\newline
We can associate the following three Hamiltonians to the line element (7.13) 
\begin{equation}
 Q({t})=\frac{1}{2}G^{\Lambda\Pi}P_{\Lambda}P_{\Pi},
\end{equation}
\begin{equation}
\widehat{H}=\frac{1}{2}\exp(2\sigma)\eta^{(\mathbf{A})(\mathbf{B})}P_{(\mathbf{A})}P_{(\mathbf{B})},
\end{equation}
and
\begin{equation}
{H}=\frac{1}{2}\eta^{\mathbf{A}\mathbf{B}}P_{\mathbf{A}}P_{\mathbf{B}}.
\end{equation}
Mapping among Hamiltonians of the form (7.16) is simpler than mapping among the forms (7.14) or (7.15).
\renewcommand{\theequation}{\thesection.\arabic{equation}}
\section{ Solutions of Mapping Among Manifolds } 
\setcounter{equation}{0} $         $
In this Section, we construct an exact solution of mapping among two pseudo-Riemannian manifolds.
\newline
Let us consider the pseudo-Riemannian line element (7.13),
\begin{equation}
ds^2=G_{\Lambda\Pi}du^{\Lambda}du^{\Pi}=\exp(2\sigma)\eta_{(\mathbf{A})(\mathbf{B})}dz^{(\mathbf{A})}dz^{(\mathbf{B})}=
\eta_{\mathbf{A}\mathbf{B}}dy^{\mathbf{A}}dy^{\mathbf{B}}.
\end{equation}
We can associate the following three Hamiltonians to the line element (8.1) 
\begin{equation}
 Q=\frac{1}{2}G^{\Lambda\Pi}p_{\Lambda}p_{\Pi},
\end{equation}
\begin{equation}
\widehat{H}=\frac{1}{2}\exp(2\sigma)\eta^{(\mathbf{A})(\mathbf{B})}P_{(\mathbf{A})}P_{(\mathbf{B})},
\end{equation}
and
\begin{equation}
{H}=\frac{1}{2}\eta^{\mathbf{A}\mathbf{B}}P_{\mathbf{A}}P_{\mathbf{B}}.
\end{equation}
It is easy to show that
\begin{equation}
 Q=\widehat{H}={H}.
\end{equation}
Let us consider pseudo-Riemannian manifolds with dimension $n=4$, so that, in the Hamiltonians (8.3), (8.4) and (8.5) we have indexes $(A)$ and ${\Lambda}$ running from 1 to 4 and index $A$ running from 1 to 6.
\newline
The three Hamiltonians are quadratic and we can use the results developed in section 3. However the $Z$ and ${\bar{Y}}$ matrices in (3.9) are not constants for ${Q}$ and ${\widehat{H}}$. This makes it difficult to obtain analytical and exact solutions. In practice, only numerical computation is feasible. But for ${{H}}$, given by (8.4),  $Z$ and ${\bar{Y}}$ are constant matrices and it is possible to calculate exact solutions for (8.4). We could then use these solutions together with the results developed in section 7 and obtain exact solutions of (8.2) and (8.3). In other words, exact solutions of (8.4) can be used to get exact solutions of (8.2) and (8.3). 
\newline
Let us consider two Hamiltonians on the form (8.4) as follows
\begin{equation}
{H}=\frac{1}{2}{\eta}^{\mathbf{A}\mathbf{B}}P_{\mathbf{A}}P_{\mathbf{B}}.
\end{equation}
\begin{equation}
{\bar{H}}=\frac{1}{2}{\bar{\eta}}^{\mathbf{A}\mathbf{B}}{\bar{P}}_{\mathbf{A}}{\bar{P}}_{\mathbf{B}}.
\end{equation}
In Section 2 we put $\frac{d{t}}{d\tau}Y_{ml}=\overline{Y}_{ml}$. 
\newline
Let us split the matrices into blocks.
\newline
For ${Y}$ and $Z$ we have,
\begin{equation}
\left(%
\begin{array}{cc}
  Y_1 & Y_2 \\
  Y_3 & Y_4 \\
\end{array}%
\right)
\end{equation}
\begin{equation}
\left(%
\begin{array}{cc}
  Z_1 & Z_2 \\
  Z_3 & Z_4 \\
\end{array}%
\right)
\end{equation}
For (8.7) we have ${Y_1=Y_2=Y_3=0}$,
and $Y_4=diag({\bar{\eta}}^{\mathbf{A}\mathbf{B}})$.
\newline
For (8.6) we have ${Z_1=Z_2=Z_3=0}$,
and $Z_4=diag({\eta}^{\mathbf{A}\mathbf{B}})$.
\newline
Let us rewrite (3.8) as follows,
\begin{equation}
\frac{d{{T}^i}_j}{d\tau}+{{T}^i}_k{J}^{kl}Z_{lj}=J^{im}\frac{d{t}}{d\tau}Y_{ml}{{T}^j}_k,
\end{equation}
in which we choose (3.21), $B=C=J$. Then the equations of motion will be given by,
\begin{equation}
\frac{d{{T}}_1}{d\tau}=\frac{d{t}}{d\tau}Y_{4}{{{T}}_3},
\end{equation}
\begin{equation}
\frac{d{{T}}_2}{d\tau}=\frac{d{t}}{d\tau}Y_{4}{{T}}_4-{{T}}_1Z_{4},
\end{equation}
\begin{equation}
\frac{d{{T}}_3}{d\tau}=0,
\end{equation}
and
\begin{equation}
\frac{d{{T}}_4}{d\tau}=-{{T}}_3Z_{4}.
\end{equation}
The exact solutions of the systems (811)-(8.14) are given by \begin{equation}
{T}_3=const.,
\end{equation}
\begin{equation}
{T}_1=Y_4{T}_3{t},
\end{equation}
\begin{equation}
{T}_2=-Y_4{T}_3Z_4({t\tau}),
\end{equation}
\begin{equation}
{T}_4={T}_3Z_4{\tau}.
\end{equation}
Notice that the matrix T are $(12X12)$, and the matrices $T(i)$, $Z(4)$ and $Y(4)$ are $(6X6)$.
\newline
Let us rewrite (3.4) as follows
\begin{equation}
{\eta}^A={{T}^A}_B{\xi}^B,
\end{equation}
in which
${\eta}^A=(Y^A,P_A)$, and 
\begin{equation}
{\bar{Y}}^A={{T_1}^A}_B{Y}^B+{{T_2}^A}_{(B+6)}{P}_B
\end{equation}
\begin{equation}
{\bar{P}}^A={{T_3}^{(A+6)}}_B{Y}^B+{{T_4}^{(A+6)}}_{(B+6)}{P}_B
\end{equation}
in which we have used a convenient notation for matrices elements.
\newline
Inverting (7.6) we will have for local coordinates of two pseudo-Riemannian manifolds,
\begin{equation}
{{u}^{\Lambda}=\exp(-{\sigma}}){{E}}{^{\Lambda}_{(A)}}{{Y}}^{(A)},
\end{equation}
and
\begin{equation}
{\bar{u}^{\Lambda}=\exp}(-{\bar{\sigma}}){\bar{E}}{^{\Lambda}_{(A)}}{\bar{Y}}^{(A)},
\end{equation}
in which we recall that indexes $(A)$ and ${\Lambda}$ running from 1 to 4 and index $A$ running from 1 to 6.
\newline
Let us consider a subset of (8.20) given by
\begin{equation}
{\bar{Y}}^{(A)}={{T_1}^{(A)}}_B{Y}^B+{{T_2}^{(A)}}_{(B+6)}{P}_B
\end{equation}
Substituting (8.24) into (8.23) we get,
\begin{eqnarray}
{\bar{u}^{\Lambda}}=\exp(-{\bar{\sigma}}){\bar{E}}^{\Lambda}_{(A)}[{{T_1}^{(A)}}_B{Y}^B+{{T_2}^{(A)}}_{(B+6)}{P}_B].
\end{eqnarray}
Let us define the following (4X6) matrices
\begin{equation}
W_1{^{\Lambda}}_{B}=\exp(-{\bar{\sigma}}){\bar{E}}^{\Lambda}_{(A)}{{T_1}^{(A)}}_B
\end{equation}
and
\begin{equation}
W_2{^{\Lambda}}_{(B+6)}=\exp(-{\bar{\sigma}}){\bar{E}}^{\Lambda}_{(A)}{{T_2}{^{(A)}}}_{(B+6)},
\end{equation}
\newline
then we can rewrite (8.25) as follows,
\begin{equation}
{\bar{u}^{\Lambda}}=W_1{^{\Lambda}}_{B}{Y}^{B}+W_2{^{\Lambda}}_{(B+6)}{P}_{B}.
\end{equation}
From Hamilton's equation,
\begin{equation}
{\bar{p}_{\Pi}}={\bar{G}_{\Pi\Lambda}}\frac{d{{\bar{u}^{\Lambda}}}}{dt}.
\end{equation}
Using (8.28) in (8.29) and defining the following matrices elements
\begin{equation}
M_1{^{\Lambda}}_{B}=\frac{d{[\exp(-\bar{\sigma}){\bar E}^{\Lambda}_{(A)}}]}{dt}{{T_2}^A}_B,
\end{equation}
\begin{equation}
M_3{^{\Lambda}}_{B}=\exp(-\bar{\sigma}){\bar E}^{\Lambda}_{(A)}Y_4{{T_3}^{(A+6)}}_B.
\end{equation}
\begin{equation}
M_2{^{\Lambda}}_{(B+6)}=\frac{d{[\exp(-\bar{\sigma}){\bar E}^{\Lambda}_{(A)}}]}{dt}{{T_2}^{(A)}}_{(B+6)},
\end{equation}
\begin{equation}
M_4{^{\Lambda}}_{(B+6)}=\exp(-\bar{\sigma}){\bar E}^{\Lambda}_{(A)}Y_4{{T_4}^{(A+6)}}_{(B+6)}.
\end{equation}
Let us set,
\begin{equation}
W_3{^{\Lambda}}_{B}=M_1{^{\Lambda}}_{B}+M_3{^{\Lambda}}_{B}
\end{equation}
\begin{equation}
W_4{^{\Lambda}}_{(B+6)}=M_2{^{\Lambda}}_{(B+6)}+M_4{^{\Lambda}}_{(B+6)}.
\end{equation}
Then equation (8.29) can be rewritten as
\begin{equation}
{\bar{p}_{\Pi}}={\bar{G}_{\Pi\Lambda}}[W_3{^{\Lambda}}_{B}Y^{B}+W_4{^{\Lambda}}_{(B+6)}P_{B}].
\end{equation}
$\Lambda$ and $\Pi$ are tensor indices, so we can use the metric tensor as follows
\begin{equation}
W_3{_{\Lambda}}_{B}={\bar{G}_{\Pi\Lambda}}W_3{^{\Lambda}}_{B},
\end{equation}
\begin{equation}
W_4{_{\Lambda}}_{(B+6)}={\bar{G}_{\Pi\Lambda}}W_4{^{\Lambda}}_{(B+6)},
\end{equation}
Then we can rewrite (8.36) as
\begin{equation}
{\bar{p}_{\Pi}}=W_3{_{\Lambda}}_{B}Y^{B}+W_4{_{\Lambda}}_{(B+6)}P_{B}.
\end{equation}
In order to analyze the meaning of $(8.28)$ and $(8.39)$, it is necessary to rewrite the line elements of two different manifolds, their corresponding Hamiltonians as well $(8.28)$ and $(8.39)$, as follows
\begin{eqnarray}
ds^2= G_{\Lambda\Pi}du^{\Lambda}du^{\Pi},
\end{eqnarray}
\begin{eqnarray}
ds^2= \bar{G}_{\Lambda\Pi}d\bar{u}^{\Lambda}d\bar{u}^{\Pi},
\end{eqnarray},
\begin{equation}
 Q=\frac{1}{2}G^{\Lambda\Pi}p_{\Lambda}p_{\Pi},
\end{equation}
\begin{equation}
 \bar{Q}=\frac{1}{2}\bar{G}^{\Lambda\Pi}\bar{p}_{\Lambda}\bar{p}_{\Pi},
\end{equation}
\begin{equation}
{\bar{u}^{\Lambda}}=W_1{^{\Lambda}}_{B}{Y}^{B}+W_2{^{\Lambda}}_{(B+6)}{P}_{B}.
\end{equation}
and
\begin{equation}
{\bar{p}_{\Pi}}=W_3{_{\Lambda}}_{B}Y^{B}+W_4{_{\Lambda}}_{(B+6)}P_{B}.
\end{equation}
Equations (8.44) and (8.45) are  transformations of coordinates and momenta of two phase spaces associated with two different manifolds given by the line elements (8.40) and (8.41).
\renewcommand{\theequation}{\thesection.\arabic{equation}}
\section{ Solutions of Matrices Riccati Equations }  
\setcounter{equation}{0} $         $
In this Section, we construct a set with an infinite number of exact solutions  of matrix Riccati equations.
\newline
Let us rewrite the system (8.11)-(8.14),
\begin{equation}
{T}_3=const.,
\end{equation}
\begin{equation}
{T}_1=Y_4{T}_3{t},
\end{equation}
\begin{equation}
{T}_2=-Y_4{T}_3Z_4({t\tau}),
\end{equation}
\begin{equation}
{T}_4=-{T}_3Z_4{\tau}.
\end{equation}
We remember that the matrices  $T(i)$, $Z(4)$ and $Y(4)$ are $(6X6)$.
\newline
Let us introduce the matrices S, R, D and E of the systems (3.17) and (3.18) and decompose each of them as the  matrix ${T}$ given by (3.10). For such, we identify $L_i=(S_i, R_i, D_i, E_I)$, in which we have index $i$ running from 1 to 4,
\begin{equation}
\left(%
\begin{array}{cc}
  L_1 & L_2 \\
  L_3 & L_4 \\
\end{array}%
\right)
\end{equation}
Let us introduce the solutions of the matrices S, R, D and E,
\begin{equation}
{S}_1=Y_4{S}_3{t},
\end{equation}
\begin{equation}
{S}_2=Y_4{S}_4{t},
\end{equation}
\begin{equation}
{R}_2=-{R}_1Z_4{\tau}.
\end{equation}
\begin{equation}
{R}_4=-{R}_3Z_4{\tau}.
\end{equation}
\begin{equation}
{D}_1=Y_4{D}_3{t},
\end{equation}
\begin{equation}
{D}_2=Y_4{D}_4{t},
\end{equation}
\begin{equation}
{E}_2=-{E}_1Z_4{\tau},
\end{equation}
and
\begin{equation}
{E}_4=-{E}_3Z_4{\tau}.
\end{equation}
\newpage
The system (9.6)-(9.13) is quite general, but the choice of some solutions for matrices S and R, for example, affects the solutions of matrices D and E because the systems (3.17) and (3.18).
\newline 
For a better understanding, let us suppose  systems (3.17) and (3.18) reduce into systems (2.21) and (2.22). In this case, it is easily seen that the matrices $S_3$, $S_4$, $R_1$ and $R_3$ are necessarily constant as are matrices $A$ in 
\newline
$T=SAR$.  
\newline
As we can always choose S and R as being non-singular, it will always be possible to use the system (3.22) in the form (3.24).
\newline
Substituting the system (9.6)-(9.13) in (3.22), assuming the matrices $S_3$, $S_4$, $R_1$, $R_3$ as constant and $A_1$, $A_2$, $A_3$, $A_4$ as not constant, we have,
\begin{equation}
{A}_1=-(S_3)^{-1}{S}_4A_3+(S_3)^{-1}\int_{a}^{b} f(\tau) \,d\tau\ +V_1,
\end{equation}
\begin{equation}
{A}_2=-(S_3)^{-1}{S}_4A_4+(S_3)^{-1}\int_{a}^{b} g(\tau) \,d\tau\ +V_2,
\end{equation}
in which matrices $V_1$ and $V_2$ are constant and  functions $f(\tau)$ and $g(\tau)$ are arbitrary. From these conditions, we conclude matrices $A_3$, $A_4$ are arbitrary.
\newline
Let us consider a second set of solutions for the matrices $A_i$,
\begin{equation}
{A}_1=-A_2R_3(R_1)^{-1}+(R_1)^{-1}\int_{a}^{b} l(\tau) \,d\tau\ +U_1,
\end{equation}
\begin{equation}
{A}_3=-A_4R_3(R_1)^{-1}+(R_1)^{-1}\int_{a}^{b} m(\tau) \,d\tau\ +U_2,
\end{equation}
in which matrices ${U_1}$ and ${U_2}$ are constant and  functions ${l}$, ${m}$ are arbitrary. The same occurs with matrices ${A}_2$, ${A}_4$, they are also arbitrary.
\newline 
If we replace (9.1)-(9.4), (9.6)-(9.13) and (9.14)-(9.15) (or (9.16)-(9.17)) on $T=SAR$, we get an equality. 
\newline
The choice $S_3$, $S_4$, $R_1$ and $R_3$ constants was arbitrary and allowed the calculation of exact solutions of some matrix equations. From the point of view of numerical computation, it will be possible a set of solutions much larger than we have obtained.
\newpage
The choice $S_3$, $S_4$, $R_1$ and $R_3$ constants will simplify the D and E matrices as follows,
\begin{equation}
{D}_3=-[S_3(A_1F_1+A_2F_3)+S_4(A_3F_1+A_4F_3)],
\end{equation}
\begin{equation}
{D}_4=-[S_3(A_1F_2+A_2F_4)+S_4(A_3F_2+A_4F_4)],
\end{equation}
\begin{equation}
{E}_1=-[G_1(A_1R_1+A_2R_3)+G_2(A_3R_1+A_4R_3)],
\end{equation}
\begin{equation}
{E}_3=-[G_3(A_1R_1+A_2R_3)+G_4(A_3R_1+A_4R_3)].
\end{equation}
in which $F_i$ and G$_i$ are arbitrary matrices and $A_i$ are given by (9.14)-(9.15) or (9.16)-(9.17).
\newpage
\renewcommand{\theequation}{\thesection.\arabic{equation}}
\section{ Concluding Remarks } 
One objective of this paper was to obtain information about a system of differential equations in which the solutions are unknown. For such, we use another system whose solutions are known. However, W matrices are functions of coordinates and momenta of the two systems, one unknown and the other known. This restricts our initial objective, but there is a large number of systems of differential equations in which only the coefficients of derivatives are known and this might be a good opportunity to use this formalism. On the other hand, for well-known situations, such as Schwarzschild and Reissner-Nordstron metrics, we will be able to construct the mapping between these important geometries.
\newline
Another important objective of this paper was to offer a set with an infinite number of exact solutions of the Riccati quadratic matrix equation. It was also shown that this could be associated with a set with an infinite number of Hamiltonians. It is always possible to associate matrices Riccati equations with a set of infinite number of Hamiltonians.


\newpage
\begin{thebibliography}{99}
\bibitem{1} A. C.V.V.de Siqueira, Matrix Riccati Equations, Kaluza-Klein, Finsler Spaces, and Mapping Among Manifolds, Preprint March 2021, {\it ResearchGate.}
\bibitem{2} P. G. L. Leach, On the theory of time-dependent linear canonical transformations as applied to Hamiltonians of the harmonic oscillator type,  {\it J. Math. Phys.} {18}, 1608 (1977).
\bibitem{3}E. A. Coddington and N. Levinson, {\bf Theory of Ordinary Differential Equations}. McGraw-Hill, New York, 1955.
\bibitem{4} E. Calabi, Improper affine hyperspheres of convex type and a generalization of a theorem by K. Jogens. Michigan Math. J. 5, 2 (1958).
\bibitem{5} A. V. Pogorelov, {\bf The Minkowski Multidimensional Problem }, John Wiley and Sons, 1978.
\end{thebibliography}
\end{document}